\documentclass{agtart_a}
\pdfoutput=1

\title{A geometric proof that $\mathrm{SL}_2(\mathbb{Z}[t,t^{-1}])$ is not
finitely presented}

\author{Kai-Uwe Bux}
\givenname{Kai-Uwe}
\surname{Bux}
\address{Department of Mathematics\\
University of Virginia\\\newline
PO Box 400137\\
Charlottesville VA 22094-4137\\
USA}
\urladdr{}

\author{Kevin Wortman}
\givenname{Kevin}
\surname{Wortman}
\address{Mathematics Department\\
Yale University\\\newline
PO Box 208283\\
New Haven CT 06520-8283\\
USA}
\urladdr{}

\volumenumber{6}
\issuenumber{}
\publicationyear{2006}
\papernumber{31}
\startpage{839}
\endpage{852}

\doi{}
\MR{}
\Zbl{}


\received{15 December 2004}
\revised{11 April 2006}
\accepted{28 October 2005}
\published{11 July 2006}
\publishedonline{11 July 2006}
\proposed{}
\seconded{}
\corresponding{}
\editor{}
\version{}

\arxivreference{math.GR/0412101}

\input basic_notation

\let\xysavmatrix\xymatrix
\def\xymatrix{\disablesubscriptcorrection\xysavmatrix}
\AtBeginDocument{\let\bar\wbar\let\tilde\wtilde\let\hat\what}

\newcommand{\ZZZ}{\mathbb{Z}}
\newcommand{\FFF}{\mathbb{F}}
\newcommand{\PPP}{\mathbb{P}}
\newcommand{\CCC}{\mathbb{C}}
\newcommand{\NNN}{\mathbb{N}}
\newcommand{\QQQ}{\mathbb{Q}}
\newcommand{\RRR}{\mathbb{R}}
\newcommand{\Oka}{\mathcal{O}}
\newcommand{\notion}{\emph}

\newcommand{\Lquot}[2]{%
\begingroup
  \setbox254=\vbox{$#1$}
  \setbox253=\vbox{$#2$}
  \dimen255=\ht254
  \advance\dimen255 by \ht253%
  \dimen254=\dimen255
  \divide\dimen254 by 2%
  \dimen253=1ex
  \advance\dimen253 by \dimen254
  \advance\dimen253 by -\ht254%
  \dimen252=\dimen253
  \advance\dimen252 by -\ht253%
  \dimen251=\dimen253
  \advance\dimen251 by \dimen252
  \divide\dimen251 by 2%
  \mbox{%
    \(
      {\raisebox{\dimen252}{$#1$}}%
      \mkern-3mu%
      \raisebox{\dimen251}{$\left\backslash\rule{0pt}{0.75\dimen255}\right.$}%
      \mkern-7mu%
      {\raisebox{\dimen253}{$#2$}}
    \)%
  }%
\endgroup
}

\newcommand{\Datum}[3]{%
  \begingroup
    \year=#1%
    \month=#2%
    \day=#3%
    {\today}%
  \endgroup
}
\def\@datum#1/#2/#3.{%
  \Datum{#1}{#2}{#3}%
}
\newcommand{\datum}[1]{%
  \expandafter\@datum#1.%
}

\makeatletter
\def\cnewtheorem#1[#2]#3{\newtheorem{#1}{#3}[section]
\expandafter\let\csname c@#1\endcsname\c@lemma}
\def\dnewtheorem#1[#2]#3{\newtheorem{#1}{#3}
\expandafter\let\csname c@#1\endcsname\c@qprop}

\newtheorem{lemma}{Lemma}[section]
\cnewtheorem{observation}[lemma]{Observation}
\cnewtheorem{cor}[lemma]{Corollary}
\cnewtheorem{prop}[lemma]{Proposition}
\newtheorem{qprop}{Proposition}
\dnewtheorem{qques}[qprop]{Question}

\makeatother

\makeautorefname{observation}{Observation}
\makeautorefname{qprop}{Proposition}
\makeautorefname{qques}{Question}

\newenvironment{enumth}{\enumerate}{\endenumerate}

\newcommand{\tlquot}[2]{#2/#1}
\newcommand{\pref}[1]{\eqref{#1}}

\newcommand{\ourparagraph}[1]{\medskip\par\textbf{#1}\qua}

\newvariable{\TheMatrixSize}{n}
\newvariable{\Sl}{\operatorname{SL}}
\newvariable{\Gl}{\operatorname{GL}}
\newcommand{\subring}{\leq}

\newvariable{\GaloisField}{\FFF}
\newvariable{\ThePrimePower}{q}

\newvariable{\TheCurve}{C}
\newvariable{\TheField}{k}
\newcommand{\RatFunctions}{\TheFieldOf{\TheCurve}}

\newcommand{\OurRingOf}[1]{%
  \begingroup
    \def\tempa{#1}%
    \def\usual{\TheRoot,\AltRoot}%
    \ifx\tempa\usual
      \def\next{\IntLaurent}%
    \else
      \def\next{!!! FIXME !!!}%
    \fi
  \expandafter\endgroup\next
}
\newcommand{\OurFieldOf}[1]{%
  \begingroup
    \def\tempa{#1}%
    \def\usual{\TheRoot,\AltRoot}%
    \def\next{!!! FIXME !!!}%
    \ifx\tempa\usual
      \def\next{\RatLaurent}%
    \fi
    \def\usual{\AltRoot}
    \ifx\tempa\usual
      \def\next{\RatInvPoly}%
    \fi
    \def\usual{\TheRoot}
    \ifx\tempa\usual
      \def\next{\RatPoly}%
    \fi
  \expandafter\endgroup\next
}

\newcommand{\EilenbergMacLaneOf}[1]{\operatorname{K}(#1,1)}

\newvariable{\AlgRat}{\overline{\QQQ}}
\newvariable{\AlgInt}{\overline{\ZZZ}}
\newvariable{\TheRoot}{\infty}
\newvariable{\AltRoot}{\Zero}

\newvariable{\TheCurvePoint}{x}
\newvariable{\AltCurvePoint}{y}
\newvariable{\TheFunction}{f}

\newvariable{\TheIsometry}{l}
\newvariable{\TheReal}{s}
\newvariable{\AltReal}{s'}

\newvariable{\TheIntDomain}{J}

\newvariable{\TheNumberField}{\QQQ}
\newvariable{\TheNumberRing}{\TheIntegers}
\newvariable{\ValuationRing}{\Oka}
\newvariable{\IntRing}{\Oka}
\newvariable{\TheReals}{\RRR}
\newvariable{\TheIntegers}{\ZZZ}
\newvariable{\TheComplexNumbers}{\CCC}

\newvariable{\TheNumberOfPlaces}{m}
\newvariable{\TheIndex}{i}

\newvariable{\PosCharacteristic}{p}

\newvariable{\TheUnity}{1}

\newcommand{\Pair}[2]{(#1,#2)}

\newvariable{\Link}{\operatorname{Link}}

\newvariable{\TheBasis}{\mathcal{B}}
\newvariable{\SubBasis}{\mathcal{F}}
\newvariable{\SubSpace}{W}
\newvariable{\AffGroup}{\mathbf{G}}

\newvariable{\TheProjective}{P}
\newvariable{\AbstractGroup}{\Gamma}
\newvariable{\TheGenerator}{\xi}

\newvariable{\TheDim}{d}
\newvariable{\TheCell}{\sigma}

\newvariable{\TheUnit}{\TheVariable}
\newvariable{\TheUniformizer}{\pi}

\newvariable{\TheNumber}{n}
\newvariable{\TheInteger}{m}
\newvariable{\Homology}{\operatorname{H}}
\newvariable{\RedHomology}{\widetilde{\Homology}}

\newvariable{\TheRange}{N}
\newvariable{\TheCoefficient}{q}

\newvariable{\TheDegreeBound}{d}

\newcommand{\TheShift}[1][]{%
  \begingroup
    \def\tempa{#1}%
    \def\minus{-}%
    \ifx\tempa\minus
      \def\next{}%
    \else
      \def\next{-}%
    \fi
  \expandafter\endgroup\next
}
\newcommand{\AltShift}[1][]{
  \begingroup
    \def\tempa{#1}%
    \def\minus{-}%
    \ifx\tempa\minus
      \def\next{-}%
    \else
      \def\next{}%
    \fi
  \expandafter\endgroup\next
}
\newvariable{\DiagMatrix}{D}
\newvariable{\UnipMatrix}{U}
\newvariable{\MatrixExponent}{n}
\newvariable{\MatrixIndex}{n}
\newvariable{\MatrixPair}{A}

\newvariable{\TheSubspace}{X}
\newvariable{\AltSubspace}{Y}
\newvariable{\TheConePoint}{\mathbf{x}}
\newvariable{\ThePoint}{\mathbf{z}}
\newvariable{\TheSegment}{\sigma}
\newvariable{\TheUnipotent}{u}
\newvariable{\TheLoop}{\gamma}

\newvariable{\TheQuadrant}{Q}
\newvariable{\TheRetraction}{\rho}

\newvariable{\TheBound}{C}

\newvariable{\TheGroup}{G}

\newvariable{\TheBorel}{B}
\newvariable{\TheMonomial}{N}
\newvariable{\TheStandardReflections}{S}
\newvariable{\TheTriangle}{\Delta}
\newvariable{\TheCone}{C}

\newvariable{\TheSpace}{X}
\newvariable{\TheSpacePoint}{x}
\newvariable{\TheRadius}{r}
\newvariable{\TheBall}{\operatorname{B}}
\newvariable{\Nbhd}{\operatorname{Nbhd}}

\newvariable{\IrrPolynomial}{p}
\newvariable{\ThePolynomial}{r}
\newvariable{\AltPolynomial}{s}

\newvariable{\TheRingElement}{a}
\newvariable{\AltRingElement}{b}

\newvariable{\TheInclusion}{\iota}
\newvariable{\TheProjection}{\pi}

\newvariable{\AntiDiagonal}{V}

\newvariable{\TheLine}{L}
\newvariable{\TheX}{x}
\newvariable{\TheY}{y}

\newvariable{\TheValuation}{v}
\newvariable{\TheLaurent}{q}
\newvariable{\TheLaurentCoeff}{a}

\newcommand{\Ftype}[1]{F${}_{#1}$}
\newcommand{\FPtype}[1]{FP${}_{#1}$}

\newcommand{\apply}{\cdot}
\newcommand{\crossprod}{\times}

\newvariable{\TheComplex}{X}

\newvariable{\Ccc}{\CCC}

\newvariable{\AffBuild}{X}
\newvariable{\AffApartment}{\Sigma}
\newvariable{\TheSector}{C}
\newvariable{\TheTip}{s}
\newvariable{\TheVertex}{x}
\newvariable{\AltVertex}{y}
\newvariable{\AffVertex}{\mathbf{\TheVertex}}
\newvariable{\AltAffVertex}{\mathbf{\AltVertex}}

\newvariable{\TheVectorSpace}{V}
\newvariable{\AltVectorSpace}{W}

\newcommand{\ProjectiveSpace}[1][\One]{\PPP^{#1}}

\newvariable{\TheMatrix}{M}
\newvariable{\TheSize}{\Two}
\newvariable{\SLn}{\operatorname{Sl}_{\TheSize}}
\newvariable{\GLn}{\operatorname{Gl}_{\TheSize}}
\newvariable{\TheVariable}{t}
\newvariable{\AltVariable}{s}
\newvariable{\RatRat}{\QQQ(\TheVariable)}
\newvariable{\IntRat}{\ZZZ(\TheVariable)}
\newvariable{\RatPoly}{\QQQ[\TheVariable]}
\newvariable{\RatInvPoly}{\QQQ[{\TheVariable[][-\One]}]}
\newvariable{\IntPoly}{\ZZZ[\TheVariable]}
\newvariable{\RatLaurent}{\QQQ[\TheVariable,\TheVariable^{-\One}]}
\newvariable{\IntLaurent}{\ZZZ[\TheVariable,\TheVariable^{-\One}]}

\newvariable{\TheTree}{T}

\newvariable{\OneOne}{a}
\newvariable{\OneTwo}{b}
\newvariable{\TwoOne}{c}
\newvariable{\TwoTwo}{d}
\newvariable{\TheValuationBound}{C}
\newvariable{\Diag}{D}
\newvariable{\TheShiftConstant}{a}

\newvariable{\IntArith}{\ZZZ_{\ThePrimeSet}}
\newvariable{\RatArith}{\QQQ_{\ThePrimeSet}}

\newvariable{\Deg}{\operatorname{deg}}
\newvariable{\TheCoordinate}{a}
\newvariable{\ValIndex}{i}
\newvariable{\FiltrIndex}{j}
\newvariable{\ValIndexP}{(\ValIndex)}
\newvariable{\TheVertInd}{j}
\newvariable{\AltVertInd}{k}
\newvariable{\TheLastIndex}{s}
\newvariable{\TheExponent}{n}
\newvariable{\AltExponent}{m}
\newvariable{\TheBigNumber}{N}
\newvariable{\TheType}{\TheLastIndex}

\newvariable{\ThePrimeSet}{S}
\newvariable{\TheRing}{R}

\newvariable{\RedGroup}{\mathbf{G}}

\newvariable{\TheDegree}{m}

\newvariable{\TheStabilizer}{G}
\newvariable{\TheRank}{r}

\newvariable{\TheRay}{F}

\newvariable{\Stabilizer}{\operatorname{Stab}}

\newcommand{\edge}{-\negthinspace\negthinspace\negthinspace-
		   \negthinspace\negthinspace\negthinspace-}

\renewcommand{\setminus}{-}
\renewcommand{\TwoByTwoMatrix}[4]{%
  \mathchoice{%
    \begin{pmatrix}
      #1 & #2\\
      #3 & #4
    \end{pmatrix}
  }{%
    \left(\begin{smallmatrix}
      #1 & #2\\
      #3 & #4
    \end{smallmatrix}\right)
  }{%
    \left(\begin{smallmatrix}
      #1 & #2\\
      #3 & #4
    \end{smallmatrix}\right)
  }{%
    \left(\begin{smallmatrix}
      #1 & #2\\
      #3 & #4
    \end{smallmatrix}\right)
  }
}

\begin{document}

\begin{asciiabstract}
We give a new proof of the theorem of Krstic-McCool from
the title. Our proof has potential applications to the study of
finiteness properties of other subgroups of SL_2 resulting
from rings of functions on curves.
\end{asciiabstract}

\begin{htmlabstract}
We give a new proof of the theorem of Krsti&#x0107;&ndash;McCool from
the title. Our proof has potential applications to the study of
finiteness properties of other subgroups of SL<sub>2</sub> resulting
from rings of functions on curves.
\end{htmlabstract}

\begin{abstract}
We give a new proof of the theorem of Krsti\'{c}--McCool from
the title. Our proof has potential applications to the study of
finiteness properties of other subgroups of $\mathrm{SL}_2$ resulting
from rings of functions on curves.
\end{abstract}

\maketitle

\section{Introduction}
Our main result is a strengthening of the theorem of Krsti\'{c}--McCool
from the title.
\begin{qprop}\label{propa}
  The group $\SlOf[\Two]{\IntLaurent}$
  is not finitely presented, indeed it is not even
  of type \FPtype{\Two}.
\end{qprop}
It will be clear from our proof that $\TheIntegers$ can be replaced
in \fullref{propa} with any ring of integers in an algebraic number field.
Note that the theorem of Krsti\'{c}--McCool \cite{KMcC97} also
allows for this replacement as well as for many other generalizations of
the ring $\IntLaurent$, which include
in particular
any ring of the form
${\TheIntDomain}[\TheVariable,\TheVariable[][-\One]]$ where $\TheIntDomain$
is an integral domain.

Let us recall the definition of type \FPtype{\Two}.
\ourparagraph{Type \boldmath\FPtype{\TheType}}
  A group $\AbstractGroup$ is of \notion{type \FPtype{\TheType}}
  if $\TheIntegers$, regarded as a
  $\TheIntegers\AbstractGroup$--module via the trivial action,
  admits a partial projective resolution
  \[
    \TheProjective[\TheType]
    \rightarrow
    \TheProjective[\TheType-\One]
    \rightarrow\cdots\rightarrow
    \TheProjective[\One]
    \rightarrow
    \TheProjective[\Zero]
    \rightarrow
    \TheIntegers
    \rightarrow
    \Zero
  \]
  by finitely
  generated $\TheIntegers\AbstractGroup$--modules
  $\TheProjective[\TheIndex]$.

  Every group is of type \FPtype{\Zero}. Type \FPtype{\One} is equivalent
  to the property of finite generation. Every finitely-presented group
  is of type \FPtype{\Two}, but Bestvina--Brady showed
  the converse does not hold in general
  \cite[Example 6.3(3)]{BeBr97}.

\ourparagraph{Purpose}
  In \cite{BuWo04}, we studied finiteness properties of subgroups
  of linear reductive groups arising from rings of functions on algebraic
  curves defined over finite fields. For example, we showed that
  $\Sl_{\TheMatrixSize}\bigl(\GaloisField_{\ThePrimePower}[\TheVariable]\bigr)$
  is not of type \FPtype{\TheMatrixSize-\One} and
  $\Sl_{\TheMatrixSize}\bigl(\GaloisField_{\ThePrimePower}[\TheVariable,\TheVariable[][-\One]]\bigr)$
  is not of type \FPtype{\Two(\TheMatrixSize-\One)} where
  $\GaloisField[\ThePrimePower]$ is a finite field.

  We wrote this paper to show how the techniques in
  \cite{BuWo04} might be applied to a more general class of groups.

  In this paper we stripped down the general proof of the main result from
  \cite{BuWo04} to the special case of showing that
  $\Sl_{\Two}\bigl(\GaloisField_{\ThePrimePower}[\TheVariable,\TheVariable[][-\One]]\bigr)$
  is not of type
  \FPtype{\Two}, and then made some modest alterations until we
  arrived at the proof  of \fullref{propa} presented below.

  It seems likely that more results along these lines can be proved, but
  it is not clear to us how much the results in
  \cite{BuWo04} can be generalized.
  Below we phrase a question that seems a good
  place to start.

\ourparagraph{Rings of functions on curves}
  Let $\TheCurve$ be an irreducible smooth projective curve defined over an
  algebraically closed field $\TheField$. We let $\RatFunctions$
  be the field of rational functions defined
  on $\TheCurve$, and we denote the set of nonzero elements
  of this field by $\UnitsOf{\RatFunctions}$.

  For each point $\TheCurvePoint\in\TheCurve$, there is a
  discrete valuation
  \(
    \TheValuation[\TheCurvePoint]
    \mapcolon
    \UnitsOf{\RatFunctions}
    \rightarrow
    \TheIntegers
  \)
  that assigns to any
  nonzero function $\TheFunction$ on $\TheCurve$ its vanishing
  order at $\TheCurvePoint$. Formally, we extend
  $\TheValuation[\TheCurvePoint]$ to all of $\RatFunctions$ by
  \(
    \TheValuationOf[\TheCurvePoint]{\Zero} := \infty.
  \)

  We let
  $\ThePrimeSet[\One],\ThePrimeSet[\Two],\ldots,
  \ThePrimeSet[\TheNumberOfPlaces]\subseteq\TheCurve$
  be collections of pairwise disjoint
  finite nonempty sets of closed points in $\TheCurve$. We call a
  ring $\TheRing\subring\RatFunctions$ containing
  some nonconstant function and the constant function $\One$
  an \notion{$\TheNumberOfPlaces$--place ring}
  if the following two conditions are satisfied:
  \begin{enumth}
    \item\label{cond:smooth}
      For all $\TheFunction\in\TheRing$ and all
      \(
	\TheCurvePoint\in\TheCurve\setminus\Parentheses{
	  \Union[\TheIndex=\One][\TheNumberOfPlaces]{
	    \ThePrimeSet[\TheIndex]
	  }
	}
	,
      \)
      we have
      $\TheValuationOf[\TheCurvePoint]{\TheFunction} \geq \Zero$.
    \item\label{cond:poles}
      If there is an $\TheIndex$, an
      $\TheCurvePoint\in\ThePrimeSet[\TheIndex]$, and an
      $\TheFunction\in\TheRing$ such that
      \(
	\TheValuationOf[\TheCurvePoint]{\TheFunction} < \Zero
	,
      \)
      then
      \(
	\TheValuationOf[\AltCurvePoint]{\TheFunction} < \Zero
      \)
      for all $\AltCurvePoint\in\ThePrimeSet[\TheIndex]$.
  \end{enumth}

  For example,
  if $\ProjectiveSpace[\One]$ is the projective line, then
  $\TheFieldOf{\ProjectiveSpace[\One]}$ is
  isomorphic to the field
  \(
    \TheFieldOf{\TheVariable}
  \)
  of rational functions in one variable.
  Thus, if $\TheIntDomain$ is a subring of $\TheField$, then
  \(
    \TheIntDomainAd{\TheVariable}
    \subring\TheFieldOf{\ProjectiveSpace[\One]}
  \)
  is a $\One$--place ring
  with $\ThePrimeSet[\One]=\SetOf{\infty}$, while
  $J[\TheVariable,\TheVariable[][-\One]]$
  is a $\Two$--place ring with
  $\ThePrimeSet[\One]=\SetOf{\infty}$ and
  $\ThePrimeSet[\Two]=\SetOf{\Zero}$.
  For an example of a $\One$--place ring $\TheRing$ that obeys
  condition~2 nontrivially, we can take
  \(
    \TheRing = \mathbb{Z}\bigl[\frac{\One}{\TheVariable[][\Two]-\Two}\bigr]
    \subring
    \CccOf{\TheVariable}
  \)
  with
  \(
    \ThePrimeSet[\One]
    =
    \bigl\{\sqrt{\Two},-\sqrt{\Two}\bigr\}.
  \)

  Note that the definition of an $\TheNumberOfPlaces$--place ring is a
  generalization of the definition of a ring of $\ThePrimeSet$--integers
  of a global function field.

\ourparagraph{Finiteness properties of linear groups}
  We ask the following question:
  \begin{qques}\label{quesb}
    Is there an example of an $\TheNumberOfPlaces$--place ring $\TheRing$
    such that $\SlOf[\TheMatrixSize]{\TheRing}$ is of type
    \FPtype{\TheNumberOfPlaces\Parentheses{\TheMatrixSize-\One}}\,?
  \end{qques}
  Specifically,
  is there an $\TheMatrixSize\geq \Two$ such that
  $\SlOf[\TheMatrixSize]{\IntPoly}$ is of type \FPtype{\TheMatrixSize-\One}
  or such that
  $\SlOf[\TheMatrixSize]{\IntLaurent}$ is of type
  \FPtype{\Two\Parentheses{\TheMatrixSize-\One}}\,?

  There seems to be no known example as above, though
  relatively few candidates have been examined
  for this property.
  Krsti\'{c}--McCool
  \cite{KMcC97,KMcC99} proved that
  \(
    \SlOf[\Two]{J[\TheVariable,\TheVariable[][-\One]]}
  \)
  and
  \(
    \SlOf[\Three]{
      \TheIntDomainAd{
	\TheVariable
      }
    }
  \)
  are not finitely presented for any
  integral domain $\TheIntDomain$.
  In \cite{BuWo04}, we prove that there exist no
  examples when $\TheRing$ is a ring of $\ThePrimeSet$--integers
  of a global function field. Examples of such rings include
  \(
    \GaloisFieldAd[\ThePrimePower]{\TheVariable}
  \)
  and
  \(
    \GaloisField_{\ThePrimePower}[\TheVariable,\TheVariable[][-\One]].
  \)

  We also  know that there are no examples as asked for in
  \fullref{quesb} when $\TheNumberOfPlaces=\One$ and $\TheMatrixSize=\Two$.
  We give a proof of this fact in
  \fullref{sec:finite_generation}.
  This is an easy result, but as this general problem has not been studied
  extensively, it appears not to have been stated in this form
  in the literature.

\ourparagraph{About the proof}
  Our proof of \fullref{propa} is geometric in that it employs the action of
  $\SlOf[\Two]{\IntLaurent}$ on a product of
  two Bruhat--Tits trees. It is essentially a special case
  of our proof that arithmetic subgroups of
  $\Sl[\TheMatrixSize]$ over global function fields are not of
  type \FPtype{\infty} \cite{BuWo04}. The
  proof uses a result of K.\,Brown's which requires the action
  to
  have ``nice'' stabilizers.
  Unfortunately, the stabilizer types of
  \(
    \SlOf[\Two]{\TheRing}
  \)
  are unknown to us for many
  of the more interesting $\Two$--place rings $\TheRing$.
  This prevents us from applying
  our proof to groups other than
  $\Sl_{\Two}\bigl({\IntRing}[\TheVariable,\TheVariable[][-\One]]\bigr)$
  where $\IntRing$ is the ring of integers in an algebraic number field.

\ourparagraph{Other finiteness properties}
  As an aside, we point out a few loosely related facts. In
  \cite{KMcC99},
  Krsti\'{c}--McCool showed that
  \(
    \SlOf[\Three]{\TheIntDomainAd{\TheVariable}}
  \)
  is not finitely presented for any integral domain $\TheIntDomain$.
  Suslin proved in
  \cite{Susl77}
  that
  \(
    \SlOf[\TheMatrixSize]{\IntPoly}
  \)
  and
  \(
    \SlOf[\TheMatrixSize]{\IntLaurent}
  \)
  are finitely generated by elementary
  matrices when $\TheMatrixSize\geq\Three$.
  It is not known whether
  \(
    \SlOf[\Two]{\IntLaurent}
  \)
  is also
  generated by elementary matrices. In fact, even finite
  generation is an open problem for this group.

\ourparagraph{Homology}
  Our proof of \fullref{propa} can be seen as a variant of
  Stuhler's proof \cite{Stuh80} that
  \(
    \SlOf[\Two]{\GaloisField_{\ThePrimePower}
      [\TheVariable,\TheVariable[][-\One]]}
  \)
  is not of type \FPtype{\Two}.
  As Stuhler's proof establishes the stronger fact that
  the second homology
  \(
    \HomologyOf[\Two]{
      \SlOf[\Two]{\GaloisField_{\ThePrimePower}[\TheVariable,\TheVariable[][-\One]]};\TheIntegers}
  \)
  is infinitely generated, it is natural to
  wonder if the proof of \fullref{propa} below
  can be extended to show that
  \(
    \HomologyOf[\Two]{
      \SlOf[\Two]{\IntLaurent}
      ;
      \TheIntegers
    }
  \)
  is infinitely generated.

\ourparagraph{Type \boldmath\Ftype{\TheType}}
  We will not use type \Ftype{\TheType} in
  this paper, but as it is related to type \FPtype{\TheType},
  we recall its definition here.

  A group $\AbstractGroup$ is of \notion{type \Ftype{\TheType}} if
  there exists an Eilenberg--Mac\,Lane complex
  $\EilenbergMacLaneOf{\AbstractGroup}$ with finite $\TheType$--skeleton.
  For $\TheType\geq\Two$, a group is of type \Ftype{\TheType} if and
  only if it is finitely presented and of type \FPtype{\TheType}.
  In general, type \Ftype{\TheType} is stronger than
  type \FPtype{\TheType}.

\ourparagraph{Outline of the paper}
  In \fullref{sec:geometry},
  we present the main body of the proof of \fullref{propa},
  leaving the verification that cell
  stabilizers are well-behaved for
  \fullref{sec:stabilizers}. In \fullref{sec:comments},
  we comment on \fullref{quesb}.

\ourparagraph{Acknowledgments}
  We thank Benson~Farb and Karen~Vogtmann for suggesting that we
  should explore this direction. We thank
  Roger~Alperin and Kevin~P~Knudson
  for helpful conversations. We also thank the referee for suggesting
  some improvements to the paper.

\section{The action on a product of trees}\label{sec:geometry}
Let $\TheValuation[\TheRoot]$ be the degree valuation on
$\TheNumberFieldOf{\TheVariable}$ given by
\[
  \TheValuationOf[\TheRoot]{
    \frac{
      \ThePolynomialOf{\TheVariable}
    }{
      \AltPolynomialOf{\TheVariable}
    }
  }
  =
  \DegOf{\AltPolynomialOf{\TheVariable}}
  -
  \DegOf{\ThePolynomialOf{\TheVariable}},
\]
and let $\TheValuation[\AltRoot]$ be the valuation at $\Zero$,
that is, the valuation corresponding to the irreducible polynomial
$\TheVariable\in\TheNumberFieldAd{\TheVariable}$. Thus
\[
  \TheValuationOf[\AltRoot]{
    \frac{
      \ThePolynomialOf{\TheVariable}
    }{
      \AltPolynomialOf{\TheVariable}
    }
    \TheVariable[][\TheExponent]
  }
  =
  \TheExponent
\]
if $\TheVariable$ does not divide $\ThePolynomialOf{\TheVariable}$
nor $\AltPolynomialOf{\TheVariable}$.

Let $\TheTree[\TheRoot]$ (resp.\ $\TheTree[\AltRoot]$) be the
Bruhat--Tits tree associated to $\SlOf[\Two]{\TheNumberFieldOf{\TheVariable}}$
with the valuation
$\TheValuation[\TheRoot]$ (resp.\ $\TheValuation[\AltRoot]$).
We consider these trees as
metric spaces by assigning a length of $\One$ to each edge. For
a definition as well as for many of the facts we will use in this
proof, we refer to Serre's book on trees \cite{Serr77}.

\ourparagraph{Outline}
  We put
  \[
    \AffBuild := \TheTree[\TheRoot] \crossprod \TheTree[\AltRoot],
  \]
  and we let $\SlOf[\Two]{\OurRingOf{\TheRoot,\AltRoot}}$ act diagonally
  on $\AffBuild$.

  We will begin by finding an
  $\SlOf[\Two]{\OurRingOf{\TheRoot,\AltRoot}}$--invariant
  cocompact subspace
  $\TheSubspace[\Zero]\subseteq \AffBuild$.
  Then for each $\MatrixIndex \in \NNN$, we will construct a
  $\One$--cycle $\TheLoop[\MatrixIndex]$ in
  $ \TheSubspace[\Zero]$ with the property that
  for any
  $\SlOf[\Two]{\OurRingOf{\TheRoot,\AltRoot}}$--invariant
  cocompact subspace $\AltSubspace\subseteq\AffBuild$ containing
  $\TheSubspace[\Zero]$, there exists some $\MatrixIndex\in\NNN$
  such that $\TheLoop[\MatrixIndex]$ represents a nontrivial element
  of the first homology group
  $\HomologyOf[\One]{\AltSubspace}$.

  A direct application of K. Brown's filtration criterion then shows that
  $\SlOf[\Two]{\OurRingOf{\TheRoot,\AltRoot}}$
  is not of type \FPtype{\Two} as long as the cell stabilizers of the
  $\SlOf[\Two]{\OurRingOf{\TheRoot,\AltRoot}}$--action
  on $\AffBuild$ are not of type \FPtype{\Two}. We leave the
  verification of this last fact for \fullref{sec:stabilizers}.

\ourparagraph{Finding a cocompact subspace}
  A crucial part of our construction will take place in
  a flat plane inside $\AffBuild$, which we shall describe now.

  Let $\ValuationRing[\TheRoot]\subring\TheNumberFieldOf{\TheVariable}$ be
  the valuation ring associated to $\TheValuation[\TheRoot]$,
  that is, the ring of all $\TheFunction\in\TheNumberFieldOf{\TheVariable}$
  with
  \(
    \TheValuationOf[\TheRoot]{\TheFunction}\geq\Zero.
  \)
  Let
  \(
    \TheLine[\TheRoot]
    \subseteq
    \TheTree[\TheRoot]
  \)
  be the unique bi-infinite geodesic stabilized by the diagonal subgroup of
  \(
    \SlOf[\Two]{\TheNumberFieldOf{\TheVariable}}.
  \)
  We parameterize $\TheLine[\TheRoot]$ by
  an isometry
  \(
    \TheIsometry[\TheRoot] \mapcolon \RRR\rightarrow\TheLine[\TheRoot]
  \)
  such that $\TheIsometryOf[\TheRoot]{\Zero}$ is the unique vertex
  stabilized by $\SlOf[\Two]{\ValuationRing[\TheRoot]}$ and such that
  the end corresponding to the positive reals is fixed
  by all upper triangular matrices in
  \(
    \SlOf[\Two]{\TheNumberFieldOf{\TheVariable}}
  \).
  Analogously, we define
  \(
    \TheIsometry[\AltRoot] \mapcolon \RRR \rightarrow\TheLine[\AltRoot]
  \).
  The plane we shall consider is the product
  \[
    \TheLine[\TheRoot]\crossprod\TheLine[\AltRoot].
  \]

  We define a diagonal
  matrix $\DiagMatrix\in
  \SlOf[\Two]{\OurRingOf{\TheRoot,\AltRoot}}$ by:
  \[
    \DiagMatrix
    :=
    \TwoByTwoMatrix{\TheUnit}{\Zero}{\Zero}{\TheUnit[][-\One]}.
  \]
  Note that for any $\MatrixExponent\in\TheIntegers$,
  we have
  \[
    \DiagMatrix[][\MatrixExponent]
    \apply
    \Pair{
      \TheIsometryOf[\TheRoot]{\Zero}
    }{
      \TheIsometryOf[\AltRoot]{\Zero}
    }
    =
    \Pair{
      \TheIsometryOf[\TheRoot]{
	\TheShift[-]\Two\MatrixExponent
      }
    }{
      \TheIsometryOf[\AltRoot]{
	\AltShift[-]\Two\MatrixExponent
      }
    }.
  \]
  Hence, if we denote by $\AntiDiagonal$ the line
  in $\TheLine[\TheRoot]\crossprod\TheLine[\AltRoot]$ of the
  form
  \(
    \SetFamOf[
      \TheReal\in\RRR
    ]{
      \Pair{
	\TheIsometryOf[\TheRoot]{
	  \TheShift[-]\TheReal
	}
      }{
	\TheIsometryOf[\AltRoot]{
	  \AltShift[-]\TheReal
	}
      }
    },
  \)
  then $\AntiDiagonal$ has a compact image under the quotient map
  \[
    \TheProjection \mapcolon
    \AffBuild
    \longrightarrow
    \AffBuild / \Sl_{\Two}\bigl(\OurRingOf{\TheRoot,\AltRoot}\bigr).
  \]
  Note that
  \[
    \TheSubspace[\Zero]
    :=
    \TheProjectionOf[][-\One]{
      \TheProjectionOf{
	\AntiDiagonal
      }
    }
    \subseteq
    \AffBuild
  \]
  is an $\SlOf[\Two]{\OurRingOf{\TheRoot,\AltRoot}}$--invariant
  cocompact subspace.

\ourparagraph{A family of loops in \boldmath$\TheSubspace[\Zero]$}
  For any $\MatrixIndex\in\TheIntegers$,
  we define the unipotent matrix
  \(
    \UnipMatrix[\MatrixIndex]
    \in
    \SlOf[\Two]{
      \OurRingOf{\TheRoot,\AltRoot}
    }
  \)
  as
  \[
    \UnipMatrix[\MatrixIndex]
    =
    \TwoByTwoMatrix{\One}{\TheUnit[][\MatrixIndex]}{\Zero}{\One}.
  \]
  Note that $\UnipMatrix[\MatrixIndex]$ fixes a point of the form
  \(
    \Pair{
      \TheIsometryOf[\TheRoot]{\TheReal}
    }{
      \TheIsometryOf[\AltRoot]{\AltReal}
    }
    \in
    \TheLine[\TheRoot]\crossprod\TheLine[\AltRoot]
  \)
  if and only if
  \(
    \TheReal \geq \TheShift[-]\MatrixIndex
  \)
  and
  \(
    \AltReal\geq \AltShift[-]\MatrixIndex.
  \)
  Moreover, any points in the plane
  \(
    \TheLine[\TheRoot]\crossprod\TheLine[\AltRoot]
  \)
  that are not fixed
  by $\UnipMatrix[\MatrixIndex]$
  are moved outside of
  $\TheLine[\TheRoot]\crossprod\TheLine[\AltRoot]$.

  For all $\MatrixIndex\in\NNN$, we define the geodesic segment
  \(
    \TheSegment[\MatrixIndex]
    \subseteq
    \AntiDiagonal
  \)
  to be the segment with endpoints
  \(
    \Pair{
      \TheIsometryOf[\TheRoot]{
	\TheShift\MatrixIndex
      }
    }{
      \TheIsometryOf[\AltRoot]{
	\AltShift\MatrixIndex
      }
    }
  \)
  and
  \(
    \Pair{
      \TheIsometryOf[\TheRoot]{
	\TheShift[-]\MatrixIndex
      }
    }{
      \TheIsometryOf[\AltRoot]{
	\AltShift[-]\MatrixIndex
      }
    }.
  \)
  Note that $\UnipMatrix[\MatrixIndex]$ fixes the endpoint of
  $\TheSegment[\MatrixIndex]$ given by
  \(
    \Pair{
      \TheIsometryOf[\TheRoot]{
	\TheShift[-]\MatrixIndex
      }
    }{
      \TheIsometryOf[\AltRoot]{
	\AltShift[-]\MatrixIndex
      }
    }
  \)
  whereas $\UnipMatrix[-\MatrixIndex]$ fixes
  its other endpoint
  \(
    \Pair{
      \TheIsometryOf[\TheRoot]{
	\TheShift\MatrixIndex
      }
    }{
      \TheIsometryOf[\AltRoot]{
	\AltShift\MatrixIndex
      }
    }.
  \)
  Since $\UnipMatrix[\MatrixIndex]$ and
  $\UnipMatrix[-\MatrixIndex]$ commute, the union of geodesic segments
  \[
    \TheLoop[\MatrixIndex]
    :=
    \TheSegment[\MatrixIndex]
    \union
    \Parentheses{
      \UnipMatrix[\MatrixIndex]
      \apply
      \TheSegment[\MatrixIndex]
    }
    \union
    \Parentheses{
      \UnipMatrix[-\MatrixIndex]
      \apply
      \TheSegment[\MatrixIndex]
    }
    \union
    \Parentheses{
	\UnipMatrix[\MatrixIndex]
	\UnipMatrix[-\MatrixIndex]
      \apply
      \TheSegment[\MatrixIndex]
    }
  \]
  is a loop. Note that
  \(
    \TheLoop[\MatrixIndex]
    \subseteq
    \TheSubspace[\Zero].
  \)

\ourparagraph{How the loops can be filled}
  It is easy to describe a filling disc for $\TheLoop[\MatrixIndex]$
  in $\AffBuild$. Just let $\TheTriangle[\MatrixIndex]$ be the closed
  triangle with geodesic sides and vertices at the endpoints of
  $\TheSegment[\MatrixIndex]$ and at the point
  \(
    \Pair{
      \TheIsometryOf[\TheRoot]{
	\TheShift[-]\MatrixIndex
      }
    }{
      \TheIsometryOf[\AltRoot]{
	\AltShift\MatrixIndex
      }
    },
  \)
  which
  is fixed by both $\UnipMatrix[\MatrixIndex]$ and
  $\UnipMatrix[-\MatrixIndex]$. Then we define
  $\TheCone[\MatrixIndex]$ to be the union of triangles
  \[
    \TheCone[\MatrixIndex]
    :=
    \TheTriangle[\MatrixIndex]
    \union
    \Parentheses{
      \UnipMatrix[\MatrixIndex]
      \apply
      \TheTriangle[\MatrixIndex]
    }
    \union
    \Parentheses{
      \UnipMatrix[-\MatrixIndex]
      \apply
      \TheTriangle[\MatrixIndex]
    }
    \union
    \Parentheses{
	\UnipMatrix[\MatrixIndex]
	\UnipMatrix[-\MatrixIndex]
      \apply
      \TheTriangle[\MatrixIndex]
    }.
  \]
  Since $\AffBuild$ is a $\Two$--complex, it does not allow for
  simplicial $\Three$--chains (using any appropriate simplicial decomposition
  of $\AffBuild$). Since $\AffBuild$ is contractible,
  it follows that there are no nontrivial simplicial $\Two$--cycles. Hence,
  there is a unique $\Two$--chain bounding $\TheLoop[\MatrixIndex]$,
  and this consists of the simplices forming $\TheCone[\MatrixIndex]$.
  Since
  \(
    \Pair{
      \TheIsometryOf[\TheRoot]{
	\TheShift[-]\MatrixIndex
      }
    }{
      \TheIsometryOf[\AltRoot]{
	\AltShift\MatrixIndex
      }
    }
    \in
    \TheCone[\MatrixIndex],
  \)
  we have:
  \begin{lemma}\label{essential_loop}
    Each loop $\TheLoop[\MatrixIndex]\subseteq\TheSubspace[\Zero]$
    represents a nontrivial class in the first homology group of
    \(
      \AffBuild
      \setminus
      \SetOf{
	\Pair{
	  \TheIsometryOf[\TheRoot]{
	    \TheShift[-]\MatrixIndex
	  }
	}{
	  \TheIsometryOf[\AltRoot]{
	    \AltShift\MatrixIndex
	  }
	}
      }.
    \)
  \end{lemma}
    Note how our proof relies on the commutator relations
    \(
      \UnipMatrix[\MatrixIndex] \UnipMatrix[-\MatrixIndex]
      =
      \UnipMatrix[-\MatrixIndex] \UnipMatrix[\MatrixIndex]
    \)
    that were also essential in the argument of Krsti\'{c}--McCool
    \cite{KMcC97}.

\ourparagraph{An unbounded sequence in the quotient}
    We will need to know that the points
    \(
	\Pair{
	  \TheIsometryOf[\TheRoot]{\TheShift[-]\MatrixIndex}
	}{
	  \TheIsometryOf[\AltRoot]{\AltShift\MatrixIndex}
	}
    \)
    move farther and farther away from $\TheSubspace[\Zero]$.
    We will use this to show that for any
    $\SlOf[\Two]{\OurRingOf{\TheRoot,\AltRoot}}$--invariant cocompact
    subspace $\AltSubspace \subseteq \AffBuild$ containing
    $\TheSubspace[\Zero]$,
    there exists some $\MatrixIndex \in \NNN$ such that
    $\TheLoop[\MatrixIndex]$ represents a nontrivial element of
    the first homology
    group $\HomologyOf[\One]{\AltSubspace}$.

    Actually, it suffices to prove our claim for ``half of the points'':
  \begin{lemma}\label{escape}
      The sequence
      \(
	\SetFamOf[\MatrixIndex\in\NNN]{
	  \TheProjectionOf{
	    \Pair{
	      \TheIsometryOf[\TheRoot]{\TheShift[-]\Two\MatrixIndex}
	    }{
	      \TheIsometryOf[\AltRoot]{\AltShift\Two\MatrixIndex}
	    }
	  }
	}
      \)
      is unbounded in the quotient space
      \(
	  \AffBuild/\SlOf[\Two]{\OurRingOf{\TheRoot,\AltRoot}}.
      \)
  \end{lemma}
  \begin{proof}
      Note that
      \(
	\SlOf[\Two]{
	  \TheNumberFieldOf{\TheVariable}
	}
	\crossprod
	\SlOf[\Two]{
	  \TheNumberFieldOf{\TheVariable}
	}
      \)
      acts on
      \(
	\TheTree[\TheRoot]
	\crossprod
	\TheTree[\AltRoot]
      \)
      componentwise
      and recall that the valuations $\TheValuation[\TheRoot]$ and
      $\TheValuation[\AltRoot]$ define a metric on
      \(
	\SlOf[\Two]{
	  \TheNumberFieldOf{\TheVariable}
	}
	\crossprod
	\SlOf[\Two]{
	  \TheNumberFieldOf{\TheVariable}
	}
      \)
      so that vertex stabilizers are bounded subgroups.
      Thus, to prove that a set of vertices in the quotient
      \(
	  \AffBuild/\SlOf[\Two]{\OurRingOf{\TheRoot,\AltRoot}}
      \)
      is not bounded, it suffices to prove that it has an unbounded
      preimage under the canonical projection
      \[
	 \bigl(\SlOf[\Two]{
	    \TheNumberFieldOf{\TheVariable}
	  }
	  \crossprod
	  \SlOf[\Two]{
	    \TheNumberFieldOf{\TheVariable}
	  }\bigr) /
	 \Sl_{\Two}\bigl(\OurRingOf{\TheRoot,\AltRoot}\bigr)
	\longrightarrow
	  \AffBuild /
	  \Sl_{\Two}\bigl(\OurRingOf{\TheRoot,\AltRoot}\bigr)
      \]
      where
      \(
	  \SlOf[\Two]{\OurRingOf{\TheRoot,\AltRoot}}
      \)
      is embedded diagonally in
      \(
	 \SlOf[\Two]{
	    \TheNumberFieldOf{\TheVariable}
	  }
	  \crossprod
	  \SlOf[\Two]{
	    \TheNumberFieldOf{\TheVariable}
	  }.
      \)

      Put
      \(
	\MatrixPair :=
	\Pair{\DiagMatrix}{\DiagMatrix[][-\One]}
	\in
	\SlOf[\Two]{\TheNumberFieldOf{\TheVariable}}
	\crossprod
	\SlOf[\Two]{\TheNumberFieldOf{\TheVariable}},
      \)
      and observe that
      \[
	\MatrixPair[][\MatrixIndex]\apply
	\Pair{
	  \TheIsometryOf[\TheRoot]{\Zero}
	}{
	  \TheIsometryOf[\AltRoot]{\Zero}
	}
	=
	\Pair{
	  \TheIsometryOf[\TheRoot]{\TheShift[-]\Two\MatrixIndex}
	}{
	  \TheIsometryOf[\TheRoot]{\AltShift\Two\MatrixIndex}
	}.
      \]
      As we have argued, it suffices to prove that the sequence
      \(
	\SlOf[\Two]{\OurRingOf{\TheRoot,\AltRoot}}
	\MatrixPair[][\MatrixExponent]
      \)
      is unbounded in
      \(
	 \SlOf[\Two]{
	    \TheNumberFieldOf{\TheVariable}
	  }
	  \crossprod
	  \SlOf[\Two]{
	    \TheNumberFieldOf{\TheVariable}
	  }
      \)
      modulo
      \(
	  \SlOf[\Two]{\OurRingOf{\TheRoot,\AltRoot}}
	.
      \)
      So assume, for a contradiction, this sequence is bounded.
      By definition, this means that there is a global constant
      $\TheValuationBound$ satisfying the following condition:
      \begin{quote}
	For any $\MatrixExponent\in\NNN$, there is a matrix
	\(
	  \TheMatrix[\MatrixExponent]
	  =
	  \TwoByTwoMatrix%
	    {\OneOne[\MatrixExponent]}%
	    {\OneTwo[\MatrixExponent]}%
	    {\TwoOne[\MatrixExponent]}%
	    {\TwoTwo[\MatrixExponent]}%
	  \in
	  \SlOf[\Two]{\OurRingOf{\TheRoot,\AltRoot}}
	\)
	such that the values of $\TheValuation[\TheRoot]$ of
	the coefficients of
	\(
	  \TheMatrix[\MatrixExponent]
	  \DiagMatrix[][\MatrixExponent]
	\)
	are bounded from below by $\TheValuationBound$ and
	the values of $\TheValuation[\AltRoot]$ of
	the coefficients of
	\(
	  \TheMatrix[\MatrixExponent]
	  \DiagMatrix[][-\MatrixExponent]
	\)
	are also bounded from below by $\TheValuationBound$.
      \end{quote}
      Recall that
      \(
	\DiagMatrix
	=
	\TwoByTwoMatrix{\TheUnit}{\Zero}{\Zero}{\TheUnit[][-\One]}
      \)
      and that
      \(
	\TheValuationOf[\TheRoot]{\TheUnit}
	=-\One
      \)
      whereas
      \(
	\TheValuationOf[\AltRoot]{\TheUnit}
	=\One
	.
      \)
      Since
      \[
	\TheValuationBound
	\leq
	\TheValuationOf[\TheRoot]{
	  \OneOne[\MatrixExponent]
	  \TheUnit[][\MatrixExponent]
	}
	=
	\TheValuationOf[\TheRoot]{
	  \OneOne[\MatrixExponent]
	}
	+
	\MatrixExponent
	\TheValuationOf[\TheRoot]{\TheUnit}
	=
	\TheValuationOf[\TheRoot]{
	  \OneOne[\MatrixExponent]
	}
	-
	\MatrixExponent
      \]
      and
      \[
	\TheValuationBound
	\leq
	\TheValuationOf[\AltRoot]{
	  \OneOne[\MatrixExponent]
	  \TheUnit[][-\MatrixExponent]
	}
	=
	\TheValuationOf[\AltRoot]{
	  \OneOne[\MatrixExponent]
	}
	-
	\MatrixExponent
	\TheValuationOf[\AltRoot]{\TheUnit}
	=
	\TheValuationOf[\AltRoot]{
	  \OneOne[\MatrixExponent]
	}
	-
	\MatrixExponent
      \]
      we find that
      \(
	\TheValuationOf[\TheRoot]{
	  \OneOne[\MatrixExponent]
	}
	\geq\One
      \)
      and
      \(
	\TheValuationOf[\AltRoot]{
	  \OneOne[\MatrixExponent]
	}
	\geq
	\One
      \)
      whenever
      $\MatrixIndex\geq\One-\TheValuationBound$,
      which implies
      $\OneOne[\MatrixExponent]=\Zero$. However, the same argument
      shows $\TwoOne[\MatrixExponent]=\Zero$,
      for $\MatrixExponent\geq\One-\TheValuationBound$.
      But then, $\TheMatrix[\One-\TheValuationBound]\not\in
      \SlOf[\Two]{\OurRingOf{\TheRoot,\AltRoot}}$.
  \end{proof}

\ourparagraph{Brown's criterion}
  The following is an immediate consequence of
  \cite[Theorem~2.2]{Bro87a}.
  \begin{lemma}\label{Brown}
    Suppose a group $\AbstractGroup$ acts by
    cell-permuting homeomorphisms
    on a contractible CW--complex $\TheComplex$ such that
    stabilizers of $\TheDim$--cells are of type
    \FPtype{\TheType+\One-\TheDim}.
    Assume that $\TheComplex$ admits a filtration
    \[
      \TheComplex[\Zero]
      \subseteq\TheComplex[\One]
      \subseteq\TheComplex[\Two]
      \subseteq
      \cdots
      \subseteq
      \TheComplex
      =
      \Union[\FiltrIndex\in\NNN]{\TheComplex[\FiltrIndex]}
    \]
    by $\AbstractGroup$--invariant, cocompact subcomplexes
    \(
      \TheComplex[\FiltrIndex].
    \)
    Then $\AbstractGroup$ is not of
    type \FPtype{\TheType+\One} if
    each of the reduced homology homomorphisms
    \[
      \RedHomologyOf[\TheType]{\TheComplex[\Zero]}
      \longrightarrow
      \RedHomologyOf[\TheType]{\TheComplex[\FiltrIndex]}
    \]
    is nontrivial.
  \end{lemma}
  In the following section,
  we will verify that cell stabilizers of the
  $\SlOf[\Two]{\OurRingOf{\TheRoot,\AltRoot}}$--action on $\AffBuild$
  are of type \FPtype{\infty}. Assuming this hypothesis for the moment, we
  can finish the proof of \fullref{propa} as follows:
  \begin{proof}[Proof of \fullref{propa}]
    The family of loops $\TheLoop[\TheNumber]$ is contained
    within the cocompact
    subspace $\TheSubspace[\Zero]$, which is a subcomplex
    of (a suitable subdivision) of $\AffBuild$.
    Since the quotient
    \(
      \tlquot{
	\SlOf[\Two]{\OurRingOf{\TheRoot,\AltRoot}}
      }{
	\AffBuild
      }
    \)
    has countably many cells, we can extend $\TheSubspace[\Zero]$
    to a filtration
    \[
      \TheSubspace[\Zero]
      \subseteq
      \TheSubspace[\One]
      \subseteq
      \TheSubspace[\Two]
      \subseteq
      \TheSubspace[\Three]
      \subseteq
      \cdots
      \subseteq
      \AffBuild
    \]
    of $\AffBuild$ by
    $\SlOf[\Two]{\OurRingOf{\TheRoot,\AltRoot}}$--invariant,
    cocompact subcomplexes $\TheSubspace[\FiltrIndex]$.

    By \fullref{escape},
    for each index $\FiltrIndex$ there is a natural number
    $\MatrixIndex$ such that
    \[
      \TheSubspace[\FiltrIndex]
      \subseteq
      \AffBuild \setminus
      \SetOf{
	\Pair{
	  \TheIsometryOf[\TheRoot]{\TheShift[-]\MatrixIndex}
	}{
	  \TheIsometryOf[\AltRoot]{\AltShift\MatrixIndex}
	}
      }.
    \]
    Therefore, by \fullref{essential_loop}, $\TheLoop[\MatrixIndex]$
    represents a nontrivial class in
    \(
      \RedHomologyOf[\One]{\TheSubspace[\FiltrIndex]}
    \),
    thus showing that
    \[
      \RedHomologyOf[\One]{\TheSubspace[\Zero]}
      \longrightarrow
      \RedHomologyOf[\One]{\TheSubspace[\FiltrIndex]}
    \]
    is nontrivial. By
    \fullref{Brown}, $\SlOf[\Two]{\OurRingOf{\TheRoot,\AltRoot}}$ is
    not of type \FPtype{\Two}.
  \end{proof}

\section{Finiteness properties of cell stabilizers}\label{sec:stabilizers}
It remains to verify the hypothesis about cell stabilizers.
Borel and Serre \cite[11.1]{BoSe73}
have shown that arithmetic groups are of type \Ftype{\infty}.
Therefore, the following lemma proves what we need, and more:
\begin{lemma}\label{stabilizers}
  The cell stabilizers of the
  \(
    \SlOf[\Two]{\OurRingOf{\TheRoot,\AltRoot}}
  \)-action
  on $\AffBuild$ are arithmetic groups.
\end{lemma}
This section is devoted entirely to the proof of this lemma.

\begin{observation}
  The set
  \(
    \TheBasis
    :=
    \SetOf[{
      \TheUnit[][\TheExponent]
    }]{
      \TheExponent\in\TheIntegers
    }
  \)
  is a $\TheNumberField$--vector space basis for
  \(
    \OurFieldOf{\TheRoot,\AltRoot}
  \)
  such that the
  subring $\OurRingOf{\TheRoot,\AltRoot}$
  consists precisely of those elements
  in $\OurFieldOf{\TheRoot,\AltRoot}$ whose coefficients
  with respect to $\TheBasis$ are all in $\TheNumberRing$.
\end{observation}

\ourparagraph{Stabilizers of standard vertices}
  We fix the following family of \notion{standard vertices}
  in $\AffBuild$. For $\TheVertInd\in\NNN$, put
  \[
    \AffVertex[\TheVertInd]
    :=
    \Pair{
      \TheIsometryOf[\TheRoot]{\TheVertInd}
    }{
      \TheIsometryOf[\AltRoot]{\Zero}
    }.
  \]
  Recall that $\SlOf[\Two]{\TheNumberFieldOf{\TheVariable}}$ acts
  on the tree $\TheTree[\TheRoot]$.
  The vertex
  \(
    \TheIsometryOf[\TheRoot]{\TheVertInd}
    \in
    \TheTree[\TheRoot]
  \)
  has the stabilizer
  \[
    \SetOf[{
      \TwoByTwoMatrix{a}{b}{c}{d}
      \in
      \SlOf[\Two]{\TheNumberFieldOf{\TheVariable}}
    }]{
      \TheValuationOf[\TheRoot]{a},
      \TheValuationOf[\TheRoot]{d}
      \geq \Zero;\,\,
      \TheValuationOf[\TheRoot]{b}
      \geq
      -\TheVertInd;\,\,
      \TheValuationOf[\TheRoot]{c}
      \geq
      \TheVertInd
    }.
  \]
  Thus, the stabilizer
  \[
    \StabilizerOf[\OurFieldOf{\TheRoot,\AltRoot}]{
      \AffVertex[\TheVertInd]
    }
  \]
  of the vertex $\AffVertex[\TheVertInd]$
  under the diagonal
  \(
    \SlOf[\Two]{\OurFieldOf{\TheRoot,\AltRoot}}
  \)-action on
  \(
    \AffBuild
    =
    \TheTree[\TheRoot]\crossprod\TheTree[\AltRoot]
  \)
  is
  \[
    \SetOf[{
      \TwoByTwoMatrix{a}{b}{c}{d}
      \in
      \SlOf[\Two]{\OurFieldOf{\TheRoot,\AltRoot}}
    }]{
      \begin{array}{l}
	\TheValuationOf[\TheRoot]{a},
	\TheValuationOf[\TheRoot]{d}
	\geq \Zero;\,\,
	\TheValuationOf[\TheRoot]{b}
	\geq
	-\TheVertInd;\,\,
	\TheValuationOf[\TheRoot]{c}
	\geq
	\TheVertInd
	\\
	\TheValuationOf[\AltRoot]{a},
	\TheValuationOf[\AltRoot]{b},
	\TheValuationOf[\AltRoot]{c},
	\TheValuationOf[\AltRoot]{d}
	\geq
	\Zero
      \end{array}
    },
  \]
  which is an affine algebraic $\TheNumberField$--group:
  Because of the bounds on the valuations $\TheValuation[\TheRoot]$ and
  $\TheValuation[\AltRoot]$, each matrix in
  \(
    \StabilizerOf[\OurFieldOf{\TheRoot,\AltRoot}]{
      \AffVertex[\TheVertInd]
    }
  \)
  can be considered as a $\Four$--tuple
  \(
    \TupelOf{a,b,c,d}
  \)
  in the finite dimensional $\TheNumberField$--vector space
  \(
    \TheVectorSpace[\Zero]
    \crossprod
    \TheVectorSpace[\TheVertInd]
    \crossprod
    \TheVectorSpace[\overline{\TheVertInd}]
    \crossprod
    \TheVectorSpace[\Zero]
  \)
  where
    $$\TheVectorSpace[\TheVertInd]
    :=
      \SetOf[{
	\Sum[\TheIndex=\Zero][\TheVertInd]{
	  \TheCoefficient[\TheIndex]
	  \TheVariable[][\TheIndex]
	}
      }]{
	\TheCoefficient[\TheIndex]\in\TheNumberField
      },\qquad
    \TheVectorSpace[\wbar{\TheVertInd}]
    :=
      \begin{cases}
	\TheNumberField & \text{for\ } \TheVertInd = \Zero \\
	\SetOf{\Zero}   & \text{for\ } \TheVertInd > \Zero.
      \end{cases}$$
  The requirement that the determinant
  be $\One$ translates into a system of algebraic equations defining
  an affine variety in
  \(
    \TheVectorSpace[\Zero]
    \crossprod
    \TheVectorSpace[\TheVertInd]
    \crossprod
    \TheVectorSpace[\wbar{\TheVertInd}]
    \crossprod
    \TheVectorSpace[\Zero].
  \)
  This variety is an affine $\TheNumberField$--group by means of
  matrix multiplication.

  Note that the vector space $\TheVectorSpace[\TheVertInd]$ carries
  an integral structure: the lattice of integer points is
  \(
    \bigl\{
      \Sum[\TheIndex=\Zero][\TheVertInd]{
	\TheCoefficient[\TheIndex]
	\TheVariable[][\TheIndex]
      }
    ~\big|~
      \TheCoefficient[\TheIndex]
      \in
      \TheNumberRing
    \bigr\}.
  \)
  Thus, the stabilizer
  \(
    \StabilizerOf[\OurRingOf{\TheRoot,\AltRoot}]{
      \AffVertex[\TheVertInd]
    }
  \)
  of $\AffVertex[\TheVertInd]$ in
  \(
    \SlOf[\Two]{\OurRingOf{\TheRoot,\AltRoot}}
  \)
  is the arithmetic subgroup
  \[
    \SetOf[{
      \TwoByTwoMatrix{a}{b}{c}{d}
      \in
      \SlOf[\Two]{\OurRingOf{\TheRoot,\AltRoot}}
    }]{
      \begin{array}{l}
	\TheValuationOf[\TheRoot]{a},
	\TheValuationOf[\TheRoot]{d}
	\geq \Zero;\,\,
	\TheValuationOf[\TheRoot]{b}
	\geq
	-\TheVertInd;\,\,
	\TheValuationOf[\TheRoot]{c}
	\geq
	\TheVertInd
	\\
	\TheValuationOf[\AltRoot]{a},
	\TheValuationOf[\AltRoot]{b},
	\TheValuationOf[\AltRoot]{c},
	\TheValuationOf[\AltRoot]{d}
	\geq
	\Zero
      \end{array}
    }.
  \]
  The idea of the proof is to push this result forward to other
  vertices.

\ourparagraph{Other vertices are translates}
  We claim that every vertex
  \(
    \AltAffVertex
    =
    \Pair{
      \AltVertex[\TheRoot]
    }{
      \AltVertex[\AltRoot]
    }
    \in\AffBuild
  \)
  can be written as
  $\TheMatrix\apply\AffVertex[\TheVertInd]$
  for some $\TheMatrix\in\GLnOf{\OurFieldOf{\TheRoot,\AltRoot}}$
  and some $\TheVertInd\in\NNN\union\SetOf{\Zero}$.

  To see this, we will use that the ray
  \[
    \TheRay[\AltRoot]
    :=
    \TheIsometryOf[\AltRoot]{\Zero}
    \edge
    \TheIsometryOf[\AltRoot]{\One}
    \edge
    \TheIsometryOf[\AltRoot]{\Two}
    \edge
    \TheIsometryOf[\AltRoot]{\Three}
    \edge
    \cdots
  \]
  is a fundamental domain for the action of
  $\SlOf[\Two]{\OurFieldOf{\AltRoot}}$
  on $\TheTree[\AltRoot]$. This follows from the discussion in
  Serre \cite[page~86f]{Serr77} and the fact that
  $\TheVariable\mapsto\TheVariable[][-\One]$ induces a ring
  automorphism of $\RatLaurent$ that interchanges $\RatPoly$ and
  $\RatInvPoly$.

  The matrix
  \(
    \TwoByTwoMatrix{
      \TheUnit[][\AltVertInd]
    }{\Zero}{\Zero}{\One}
  \)
  translates $\TheIsometryOf[\AltRoot]{\AltVertInd}$ to
  $\TheIsometryOf[\AltRoot]{\Zero}$
  as $\TheVariable$
  is a uniformizing element for the valuation
  $\TheValuation[\AltRoot]$.
  Thus, within two moves, we can adjust the second
  coordinate of $\AltAffVertex$ to
  \(
    \TheIsometryOf[\AltRoot]{\Zero}.
  \)

  Now, we consider
  \(
    \OurFieldOf{\TheRoot}.
  \)
  In this case, the discussion in
  Serre \cite{Serr77} applies directly: the ray
  \[
    \TheRay[\TheRoot]
    :=
    \TheIsometryOf[\TheRoot]{\Zero}
    \edge
    \TheIsometryOf[\TheRoot]{\One}
    \edge
    \TheIsometryOf[\TheRoot]{\Two}
    \edge
    \TheIsometryOf[\TheRoot]{\Three}
    \edge
    \cdots
  \]
  is a fundamental domain in $\TheTree[\TheRoot]$ for the action of
  \(
    \SlOf[\Two]{\OurFieldOf{\TheRoot}}.
  \)
  This allows us to adjust the first coordinate. Note that
  every matrix in
  \(
    \SlOf[\Two]{\OurFieldOf{\TheRoot}}
  \)
  fixes the vertex
  \(
    \TheIsometryOf[\AltRoot]{\Zero}
    \in
    \TheTree[\AltRoot].
  \)
  Thus, we do not change the second coordinate during the third and
  final move.

  We conclude:
  \begin{lemma}\label{shape}
    Every vertex stabilizer in
    \(
      \SlOf[\Two]{\OurRingOf{\TheRoot,\AltRoot}}
    \)
    is of the form
    \[
      \TheMatrix
      \StabilizerOf[\OurFieldOf{\TheRoot,\AltRoot}]{\AffVertex[\TheVertInd]}
      \TheMatrix[][-\One]
      \intersect
      \SlOf[\Two]{\OurRingOf{\TheRoot,\AltRoot}}
    \]
    for some $\TheVertInd$ and some matrix
    $\TheMatrix\in\GLnOf{\OurFieldOf{\TheRoot,\AltRoot}}$.
  \end{lemma}
  We also make the following:
  \begin{observation}\label{raw_variety}
    Since multiplication by $\TheMatrix$ can lower valuations only by
    a bounded amount,
    we can find $\TheRange\in\NNN$ such that
    \[
      \TheMatrix
      \StabilizerOf[\OurFieldOf{\TheRoot,\AltRoot}]{\AffVertex[\TheVertInd]}
      \TheMatrix[][-\One]
      \subseteq
      \SetOf[{
	\TwoByTwoMatrix{a}{b}{c}{d}
	\in
	\SlOf[\Two]{\OurFieldOf{\TheRoot,\AltRoot}}
      }]{
	a,b,c,d \in \AltVectorSpace[\TheRange]
      }
    \]
    where
    \(
      \AltVectorSpace[\TheRange]
      :=
      \bigl\{
	\Sum[\TheIndex=-\TheRange][\TheRange]{
	  \TheCoefficient[\TheIndex]
	  \TheVariable[][\TheIndex]
	}
      ~\big|~
	\TheCoefficient[\TheIndex]
	\in
	\TheNumberField
      \bigr\}.
    \)
  \end{observation}

\ourparagraph{Finite dimensional approximations}
  We want to use \fullref{raw_variety} and argue that
  \(
    \TheMatrix
    \StabilizerOf[\OurFieldOf{\TheRoot,\AltRoot}]{\AffVertex[\TheVertInd]}
    \TheMatrix[][-\One]
  \)
  is an affine $\TheNumberField$--group with arithmetic subgroup
  \[
      \TheMatrix
      \StabilizerOf[\OurFieldOf{\TheRoot,\AltRoot}]{\AffVertex[\TheVertInd]}
      \TheMatrix[][-\One]
      \intersect
      \SlOf[\Two]{\OurRingOf{\TheRoot,\AltRoot}}.
  \]
  This is accomplished as follows.
  \begin{lemma}\label{variety}
    Fix $\TheRange\in\NNN$ and let
    \(
      \AffGroup
    \)
    be a $\TheNumberField$--subvariety of the affine $\TheNumberField$--variety
    \(
      \SetOf[{
	\TwoByTwoMatrix{a}{b}{c}{d}
	\in
	\SlOf[\Two]{\OurFieldOf{\TheRoot,\AltRoot}}
      }]{
	a,b,c,d \in \AltVectorSpace[\TheRange]
      }.
    \)
    Assume that $\AffGroup$ is a $\TheNumberField$--group with respect to
    matrix multiplication.
    Then
    \(
      \AffGroup
      \intersect
      \SlOf[\Two]{\OurRingOf{\TheRoot,\AltRoot}}
    \)
    is an arithmetic subgroup of $\AffGroup$.
  \end{lemma}
  \begin{proof}
    The integer points in $\AltVectorSpace[\TheRange]$ are
    \(
      \AltVectorSpace[\TheRange]
      \intersect
      \IntLaurent.
    \)
    Thus the integer points in $\AffGroup$ are
    \(
      \AffGroup\intersect
      \SlOf[\Two]{\IntLaurent}.
    \)
  \end{proof}
  We note that \fullref{variety}
  and \fullref{raw_variety} imply:
  \begin{cor}\label{arithmetic}
    All vertex stabilizers in
    \(
      \SlOf[\Two]{\OurRingOf{\TheRoot,\AltRoot}}
    \)
    are arithmetic groups.
  \end{cor}

\ourparagraph{Extending the argument to cell stabilizers}
  So far we have argued that vertex stabilizers are arithmetic.
  To extend this argument to stabilizers of cells
  of higher dimension, note
  that the action of $\SlOf[\Two]{\OurFieldOf{\TheRoot,\AltRoot}}$ on
  $\AffBuild$ is type-preserving.
  Hence the stabilizer of a cell is the
  intersection of the stabilizers of its
  vertices. To recognize such a group as
  arithmetic using the above method,
  we just have to choose $\TheRange$ large
  enough to accommodate for all the involved
  vertex stabilizers simultaneously.
  This concludes the proof of \fullref{stabilizers}.

\section[Comments on \ref{quesb}]{Comments on \fullref{quesb}}
\label{sec:finite_generation}
\label{sec:comments}
We shall
begin with answering \fullref{quesb} when $\TheNumberOfPlaces=\One$
and $\TheMatrixSize=\Two$.
\begin{prop}\label{prop:finite_generation}
  If $\TheRing$ is a $\One$--place ring, then
  \(
    \SlOf[\Two]{\TheRing}
  \)
  is not finitely generated.
\end{prop}
\begin{proof}
  By our hypothesis on $\TheRing$, there
  is an algebraically closed field $\TheField$, and an irreducible smooth
  projective curve $\TheCurve$ defined over $\TheField$ such that
  $\TheRing$ is a subring of the field
  of rational functions $\TheFieldOf{\TheCurve}$.

  Let $\ThePrimeSet[\One] \subseteq \TheCurve$
  be the finite set of closed points given in the
  definition of $\TheRing$ as a $\One$--place ring, and pick some
  $\TheCurvePoint \in \ThePrimeSet[\One]$.
  We let $\TheTree$ be the Bruhat--Tits tree for
  \(
    \SlOf[\Two]{\TheFieldOf{\TheCurve}}
  \)
  with the valuation $\TheValuation[\TheCurvePoint]$.
  We regard
  $\TheTree$ as a metric space by assigning unit length to all
  edges.

  Denote the geodesic in $\TheTree$
  corresponding to the diagonal subgroup
  of $\SlOf[\Two]{\TheFieldOf{\TheCurve}}$
  by $\TheLine$, and parameterize $\TheLine$ by an
  isometry
  \(
    \TheIsometry \mapcolon
    \TheReals \rightarrow \TheLine
  \)
  such that the end of
  $\TheLine$ corresponding to the positive reals is fixed by
  upper-triangular matrices.

  It follows from the definition of a $\One$--place ring, that there
  exists an element
  \(
    \TheFunction \in \TheRing
  \)
  such that
  \(
    \TheValuationOf[\TheCurvePoint]{\TheFunction}
    <
    \Zero.
  \)
  We use this element to
  define for each $\MatrixIndex\in\NNN$ a matrix
  \[
    \UnipMatrix[\MatrixIndex]
    :=
    \TwoByTwoMatrix{\One}{\TheFunction[][\MatrixIndex]}{\Zero}{\One}
  \]
  Note that for sufficiently large $\MatrixIndex$, there is an
  $\TheReal[\MatrixIndex]>\Zero$
  such that
  \[
    \UnipMatrix[\MatrixIndex]
    \apply
    \TheIsometryOf{
      \ClosedInterval{\Zero}{\TheReal[\MatrixIndex]}
    }
    \intersect
    \TheIsometryOf{
      \ClosedInterval{\Zero}{\TheReal[\MatrixIndex]}
    }
    =
    \SetOf{
      \TheIsometryOf{\TheReal[\MatrixIndex]}
    }.
  \]
  Note also that
  \(
    \TheReal[\MatrixIndex]
    =
    -\MatrixIndex\TheValuationOf[\TheCurvePoint]{\TheFunction}
    +
    \TheShiftConstant
  \)
  for some $\TheShiftConstant\in\TheReals$.

  We claim that for any $\TheRadius>\Zero$, the
  $\TheRadius$--metric neighborhood of the
  orbit
  \(
    \SlOf[\Two]{\TheRing} \apply \TheIsometryOf{\Zero}
    \subseteq
    \TheTree
  \)
  is not connected.
  Indeed, for large $\MatrixIndex$, the unique path
  between
  \(
    \TheIsometryOf{\Zero}
  \)
  and
  \(
    \UnipMatrix[\MatrixIndex]\apply
    \TheIsometryOf{\Zero}
  \)
  contains
  \(
    \TheIsometryOf{\TheReal[\MatrixIndex]},
  \)
  thus it suffices to show that
  \(
    \SlOf[\Two]{\TheRing} \apply
    \TheIsometryOf{\TheReal[\MatrixIndex]}
  \)
  is an unbounded sequence in
  the quotient space
  \(
    {\TheTree}/{\SlOf[\Two]{\TheRing}}.
  \)

  Observe that for each
  \(
    \MatrixIndex\in\NNN,
  \)
  the diagonal matrix
  \[
    \DiagMatrix[\MatrixIndex]
    :=
    \TwoByTwoMatrix
      {\TheFunction[][\MatrixIndex]}{\Zero}%
      {\Zero}{\TheFunction[][-\MatrixIndex]}
  \]
  acts by translations on $\TheLine$ and that
  \(
    \DiagMatrix[\MatrixIndex] \apply
    \TheIsometryOf{\Zero}
    =
    \TheIsometryOf{
      -\Two\MatrixIndex\TheValuationOf[\TheCurvePoint]{\TheFunction}
    }.
  \)
  Thus, to prove our claim it suffices to show
  that
  \(
    \SlOf[\Two]{\TheRing}
      \DiagMatrix[\MatrixIndex]\apply\TheIsometryOf{\Zero}
  \)
  is an unbounded sequence in
  \(
    \tlquot{\SlOf[\Two]{\TheRing}}{\TheTree}.
  \)

  Since point stabilizers in $\SlOf[\Two]{\TheFieldOf{\TheCurve}}$
  are bounded, we
  can further reformulate our task as showing the sequence
  $\SlOf[\Two]{\TheRing}\DiagMatrix[\MatrixIndex]$ is unbounded in
  \(
    \tlquot{\SlOf[\Two]{\TheRing}}{\SlOf[\Two]{\TheFieldOf{\TheCurve}}}.
  \)
  For this, we will employ a proof
  by contradiction: Assuming that
  \(
    \SlOf[\Two]{\TheRing}\DiagMatrix[\MatrixIndex]
  \)
  is bounded,
  there exist matrices
  \[
    \TheMatrix[\MatrixExponent]
    =
    \TwoByTwoMatrix%
      {\OneOne[\MatrixExponent]}%
      {\OneTwo[\MatrixExponent]}%
      {\TwoOne[\MatrixExponent]}%
      {\TwoTwo[\MatrixExponent]}%
    \in
    \SlOf[\Two]{\TheRing}
  \]
  such that the image of the
  matrix entries of
  \(
    \TheMatrix[\MatrixExponent]\DiagMatrix[\MatrixExponent]
  \)
  under the valuation $\TheValuation[\TheCurvePoint]$ are bounded
  from below by a constant $\TheValuationBound$. In particular,
  \[
    \TheValuationBound
    \leq
    \TheValuationOf[\TheCurvePoint]{
      \TheFunction[][\MatrixExponent]
      \OneOne[\MatrixExponent]
    }
    =
    \MatrixExponent\TheValuationOf[\TheCurvePoint]{\TheFunction}
    +\TheValuationOf[\TheCurvePoint]{\OneOne[\MatrixExponent]}.
  \]
  Since $\TheValuationOf[\TheCurvePoint]{\TheFunction}<\Zero$, it
  follows that
  \(
    \TheValuationOf[\TheCurvePoint]{\OneOne[\MatrixExponent]}
    >
    \Zero
  \)
  for all but finitely many $\MatrixExponent$.
  Combining conditions~\pref{cond:smooth} and~\pref{cond:poles}
  of the definition of a $\One$--place
  ring,
  \(
    \TheValuationOf[\AltCurvePoint]{\OneOne[\MatrixExponent]}
    \geq
    \Zero
  \)
  for all
  \(
    \AltCurvePoint\in \TheCurve.
  \)
  Therefore, $\OneOne[\MatrixExponent]$ is a
  constant function on $\TheCurve$.
  As
  \(
    \TheValuationOf[\TheCurvePoint]{\OneOne[\MatrixExponent]}
    >
    \Zero,
  \)
  we conclude that
  $\OneOne[\MatrixExponent]=\Zero$.
  Similarly, $\TwoOne[\MatrixExponent]=\Zero$
  for sufficiently large $\MatrixExponent$ which
  contradicts that $\TheMatrix[\MatrixExponent]$
  is invertible. We have completed our proof
  of the claim that for any $\TheRadius>\Zero$, the
  $\TheRadius$--metric neighborhood of
  the orbit
  \(
    \SlOf[\Two]{\TheRing} \apply \TheIsometryOf{\Zero}
    \subseteq
    \TheTree
  \)
  is not connected.

  \fullref{prop:finite_generation}
  now follows from an application of the following lemma.
\end{proof}
\begin{lemma}
  Suppose a finitely generated group $\AbstractGroup$
  acts on a geodesic metric
  space $\TheSpace$. Then, for any point $\TheSpacePoint \in \TheSpace$,
  there is a number $\TheRadius>\Zero$ such that the metric
  $\TheRadius$--neighborhood of the orbit
  of $\AbstractGroup\apply \TheSpacePoint \subseteq \TheSpace$ is connected.
\end{lemma}
\begin{proof}
  Let
  \(
    \SetOf{
      \TheGenerator[\One],\TheGenerator[\Two],\ldots,
      \TheGenerator[\TheLastIndex]
    }
  \)
  be a finite generating set for $\AbstractGroup$.
  Choose $\TheRadius$ such that the ball
  \(
    \TheBallOf[\TheRadius]{\TheSpacePoint}
  \)
  contains all translates
  \(
    \TheGenerator[\TheIndex]\apply\TheSpacePoint.
  \)
  Then
  \(
    \AbstractGroup\apply\TheBallOf[\TheRadius]{\TheSpacePoint}
    =
    \NbhdOf[\TheRadius]{
      \AbstractGroup\apply\TheSpacePoint
    }
  \)
  is connected.
\end{proof}

\ourparagraph{The question of \boldmath\FPtype{\Two}}
  After modest adjustments,
  the proofs in \fullref{sec:geometry}
  apply
  to $\SlOf[\Two]{\TheRing}$ for many other $\Two$--place
  rings $\TheRing$.
  Thus, the only obstruction to substituting one of these groups
  for $\SlOf[\Two]{\IntLaurent}$ in the proof of \fullref{propa}
  is proving results about finiteness properties of stabilizers as
  in \fullref{sec:stabilizers}.
\eject
  Certainly there are more $\Two$--place rings that produce stabilizers
  of type \FPtype{\Two} than the rings
  \(
    {\IntRing}[\TheVariable,\TheVariable[][-\One]]
  \)
  where
  $\IntRing$ is a ring of integers in an algebraic number field,
  but this is not the case for all $\Two$--place rings. For instance, this
  is clearly not the case for any uncountable ring $\TheRing$.
  For a countable example,
  consider
  \(
    {\TheIntegers}[\AltVariable,\TheVariable,\TheVariable[][-\One]]
  \)
  as the $\Two$--place ring contained in
  \(
    \overline{\TheComplexNumbersOf{\AltVariable}}(\PPP^{\One})
    \cong
    \overline{\TheComplexNumbersOf{\AltVariable}}(\TheVariable)
  \)
  where
  \(
    \overline{\TheComplexNumbersOf{\AltVariable}}
  \)
  is the algebraic closure of the field
  \(
    \TheComplexNumbersOf{\AltVariable}
  \)
  (we take
  \(
    \ThePrimeSet[\One]
    :=
    \SetOf{\Zero}
  \)
  and
  \(
    \ThePrimeSet[\Two]
    :=
    \SetOf{\infty}
  \)).
  Then the stabilizer in
  \(
    \SlOf[\Two]{
      {\TheIntegers}[\AltVariable,\TheVariable,\TheVariable[][-\One]]}
  \)
  of the ``standard vertex''
  \(
    \AffVertex[\Zero]
  \)
  in the product of Bruhat--Tits trees
  corresponding to valuations at $\Zero$ and $\infty$
  is equal to
  \(
    \SlOf[\Two]{
      \TheIntegersAd{\AltVariable}
    }
  \)
  and thus is not finitely generated
  by \fullref{prop:finite_generation}
  since
  \(
    \TheIntegersAd{\AltVariable}
  \)
  is a $\One$--place ring.

\ourparagraph{The question of higher finiteness properties}
  Note that the results of \fullref{sec:stabilizers}
  can easily be extended to
  the groups
  \(
    \SlOf[\TheMatrixSize]{\IntPoly}
  \)
  and
  \(
    \SlOf[\TheMatrixSize]{\IntLaurent}.
  \)
  Thus, the complication
  in extending our proof of \fullref{propa} to these groups lies in
  generalizing the material of \fullref{sec:geometry}.

  Of course, for the general $\TheNumberOfPlaces$--place ring $\TheRing$
  and for $\TheMatrixSize>\Two$,
  most of the details of this paper cannot be easily extended to
  \(
    \SlOf[\TheMatrixSize]{\TheRing}.
  \)

\bibliographystyle{gtart}
\bibliography{link}

\end{document}